\documentclass[12pt]{amsart}

\usepackage{graphicx, verbatim,amsmath,amssymb}

\def\vs#1 {\vskip#1truein}
\def\hs#1 {\hskip#1truein}
\newcommand{\E}{e} 
\newcommand{\EN}{E_n} 
\newcommand{\PN}{p_n}

\newcommand{\RA}{\rangle_{n,p}}

\newcommand{\be}{\begin{equation}}
\newcommand{\ee}{\end{equation}}
\newcommand{\R}{{\mathbb R}}

\newtheorem{thm}{Theorem}[section]
\newtheorem*{thm*}{Theorem}
\newtheorem{cor}[thm]{Corollary}

\everymath{\displaystyle}

\theoremstyle{definition}

\title{Surface effects in dense random graphs with sharp edge constraint}
\date{\today}

\author{
Charles Radin
\and Kui Ren 
\and Lorenzo Sadun
}
\address{Charles Radin\\Department of Mathematics\\The University of
  Texas at Austin\\ Austin, TX 78712} \email{radin@math.utexas.edu}
\address{Kui Ren\\Department of Mathematics and the Institute for Computational Engineering and Science (ICES)\\The University of
  Texas at Austin\\ Austin, TX 78712} \email{ren@math.utexas.edu}
\address{Lorenzo Sadun\\Department of Mathematics\\The University of
  Texas at Austin\\ Austin, TX 78712} \email{sadun@math.utexas.edu}

\thanks{This work was partially supported by NSF
  grants DMS-1208941, DMS-1321018, DMS-1509088 and DMS-1620473}
\begin{document}

\begin{abstract}
  We show that the random number $T_n$ of triangles in a random graph
  on $n$ vertices, with a strict constraint on the total number of
  edges, admits an expansion $T_n = an^3 + bn^2 + F_n$, where $a$ and
  $b$ are numbers, with the mean $\langle F_n \rangle = 
 O(n)$ and the standard deviation $\sigma(T_n) =\sigma(F_n)=
  O(n^{3/2})$. The presence of a `surface term' $bn^2$ has a
  significance analogous to the macroscopic surface effects of
  materials, and is missing in the model where the edge constraint 
 is removed. We also find the surface effect in other graph models using
  similar edge constraints.
\end{abstract}

\maketitle

\section{A random graph model with dependent
  edges}\label{dependentmodel}

Consider the spaces $G_n$, $n=1,\ldots,$ of simple graphs on $n$ labeled
vertices, on which we will define probability distributions,
giving us
random graph models of increasing `size' $n$. 
Let $H$ be an arbitrary but fixed graph, and for $g\in G_n$ let $T_H(g)$
denote the number of copies of $H$ found in $g$. 
We will compute the growth rates of the expectation
and variance of $T_H$, and will show that the expectation has both
`volume' and `surface' rates of growth, which are not overshadowed by
the lower rate of growth of its standard deviation. This is analogous
to the volume and surface components of macroscopic materials, and
indeed our models were chosen to mimick the statistical mechanics
model of macroscopic materials.

We define our probability distributions as follow. For each
$0<p<1$ fix some sequence $E_n = p{n\choose 2}+O(1)$ and define $\mu_n(p)$
as the uniform distribution on those $g\in G_n$ such that
the total number of edges, denoted $T_\E(g)$, is exactly equal to $E_n$.
Graphs in this ensemble are 
often called Erd\H{o}s-R\'enyi graphs \cite{ER}; see~\cite{JLR,B} for a broad overview.

Let $H$ be a fixed graph with $v$ vertices 
and $\ell>1$ edges. (For simplicity we assume every vertex in
$H$ lies on at least one edge.)
Fixing the distribution $\mu_n(p)$, we are
interested in the expectation ${\langle T_H \rangle}_{n,p}$ and the
variance $Var(T_H)_{n,p}$. 
We first need some specialized notation to simplify
the statements of the results.

Let $N_n={n \choose 2}$. For any positive integer $k$, let
$P(\EN,N_n,k) = \frac{\EN(\EN-1)\cdots(\EN+1-k)}{N_n(N_n-1)\cdots (N_n+1-k)}$.
Let $c_n$ be the number of copies of $H$ that appear in the complete
graph; specifically, $c_n=n(n-1)\cdots(n+1-v)/|S_H|$, where $S_H$ is
the group of symmetries of the graph $H$. For instance, 
if $H$ is a triangle, then
$c_n=n(n-1)(n-2)/6$, while if $H$ is a ``2-star'' 
(that is, a graph with three
vertices and two edges), then $c_n=n(n-1)(n-2)/2$.

For each integer $k$ between 0 and $\ell$, let $C_k$ be the number of
times in which two distinguishable 
copies of $H$ in the complete graph share
exactly $k$ edges. Note that \be \label{Sumk} \sum_{k=0}^\ell C_k = c_n^2;
\qquad \sum_{k=0}^\ell k C_k = \frac{\ell^2 c_n^2}{N_n}. \ee The second
equation comes from the fact that each of the $\ell$ edges in the
first copy of $H$ has probability $1/N_n$ of being the same as each of
the $\ell$ edges of the second copy.  The two sides are just different
expressions for the sum, over all configurations, of the number 
of shared edges. Note also that $C_0$ is of order $n^{2v}$, $C_1$ is
of order $n^{2v-2}$, $C_2$ and $C_3$ are of order $n^{2v-3}$, and all
other terms are of order $n^{2v-4}$ or smaller.

\begin{thm}\label{dependenttheorem} The expectation and variance of 
$T_H$ are given by 
\begin{eqnarray}\label{thm1eq}
\langle T_H \RA & = & c_n P(\EN, N_n, \ell)  \cr &=&  c_n p^\ell +
O(n^{v-2}) = \frac{1}{|S_H|}[n^v p^\ell-\frac{v(v
    -1)}{2}n^{v-1} p^\ell] +O(n^{v-2}), \cr 
Var(T_H)_{n,p} & = & C_2 p^{2v-2}(1-p)^2 + C_3 p^{2v-3}(1-3p^2+2p^3) +
O(n^{2v-4}) \cr &=&O(n^{2v-3}).
\end{eqnarray}
If $H$ does not contain any triangles, 
then the $C_3$ term in the formula for $Var(T_H)_{n,p}$ 
is itself $O(n^{2v-4})$ and can be ignored. 
\end{thm}

In particular, the standard deviation
of $T_H$ has a lower growth
rate, $O(n^{v-3/2})$, than that of the $O(n^{v-1})$ second term in the
expansion of $T_H$, implying a meaningful surface effect.

\begin{proof} The formula for the expectation is easy. 
Each of the $c_n$ configurations has probability $P(\EN,N,\ell)$ of
appearing. We also note that 
\begin{eqnarray}\label{ExpandP}
P(\EN,N_n,k) & = & \frac{\EN(\EN-1)\cdots(\EN+1-k)}{N_n(N_n-1)\cdots (N_n+1-k)} \cr 
& = & \left( \frac{\EN}{N_n}\right)^k \frac{(1 - \frac{1}{\EN}) \cdots(1 - \frac{k-1}{\EN})}{(1 - \frac{1}{N_n}) \cdots(1 - \frac{k-1}{N_n})}\cr
& = & \left( \frac{\EN}{N_n}\right)^k \frac{1 - \frac{k(k-1)}{2\EN} + O(n^{-4})}{1 - \frac{k(k-1)}{2N_n} + O(n^{-4})} \cr 
& = & \left( \frac{\EN}{N_n}\right)^k \left ( 1 + \frac{k(k-1)}{2N_n}\left (1 - \frac{N_n}{\EN} \right )\right ) + O(n^{-4}) \cr 
& = & \PN^k - \frac{k(k-1)}{2N_n} (\PN^{k-1}-\PN^k) + O(n^{-4}),
\end{eqnarray}
where $\PN\equiv E_n/{n \choose 2}$.
Since $c_n$ is $O(n^v)$, this implies that $c_n P(\EN,N_n,\ell) = c_n \PN^\ell + O(n^{v-2})$.

Next we compute the variance. The expected value of $T_H^2$ is obtained by writing down all the configurations of 
two $H$'s, and adding their probabilities. That is, 
\be \langle T_H^2 \rangle_{n,p} = \sum_{k=0}^\ell C_k P(\EN,N_n, 2 \ell - k). \ee
Meanwhile,
\begin{eqnarray} \langle T_H \rangle_{n,p}^2 & = & c_n^2 P(\EN,N_n,\ell)^2 \cr 
& = &  c_n^2 \left ( (P(\EN,N_n,\ell)^2 - P(\EN,N_n, 2\ell) \right ) + c_n^2 P(\EN,N_n,2\ell)  \cr 
& = & c_n^2  \left ( (P(\EN,N_n,\ell)^2 - P(\EN,N_n, 2\ell) \right )+ \sum_k C_k P(\EN,N_n,2\ell).
\end{eqnarray}
However, by equation (\ref{ExpandP}), 
\be P(\EN,N_n,\ell)^2 - P(\EN,N_n,2\ell) = \frac{\ell^2}{N_n} (\PN^{2\ell-1}-\PN^{2\ell}) + O(n^{-4}), \ee
but by (\ref{Sumk}), $\ell^2 c_n^2/N_n = \sum_{k=0}^\ell k C_k$, so 
\be \langle T_H \rangle_{n,p}^2 = \sum_{k=0}^\ell C_k \left [ P(\EN,N_n,2\ell) + k(\PN^{2\ell-1}-\PN^{2\ell}) + O(n^{-4}) \right ]. \ee
This makes the variance
\be Var(T_H)_{n,p} = \sum_{k=0}^\ell C_k \left [ P(\EN,N_n,2\ell-k)-P(\EN,N_n,2\ell) -k(\PN^{2\ell-1}-\PN^{2\ell})\right ] + O(n^{2v-4}). \ee
The $k=0$ term is identically zero. All the other $C_k$'s are at most $O(n^{2v-2})$, so we can use the approximations
$P(\EN,N_n,2\ell-k) =  \PN^{2\ell-k} + O(n^{-2})$ and $P(\EN,N_n,2\ell) = \PN^{2\ell}+O(n^{-2})$ to get 
\begin{eqnarray} Var(T_H)_{n,p} &=& 
\sum_{k=1}^\ell C_k \left (\PN^{2\ell-k} - \PN^{2\ell} - k(\PN^{2\ell-1}-\PN^{2\ell}) \right ) + O(n^{2v-4}) \cr 
&=& \sum_{k=1}^\ell C_k \left (p^{2\ell-k} - p^{2\ell} - k(p^{2\ell-1}-p^{2\ell}) \right ) + O(n^{2v-4}), \end{eqnarray}
since $p_n = p + O(n^{-2})$. 
In this last sum the $k=1$ term is zero, 
the $k=2$ term is $C_2 p^{2\ell-2}(1-p)^2$, the $k=3$ term is 
$C_3 p^{2\ell-3}(1-3p^2+2p^3)$, and all remaining terms 
are $O(n^{2v-4})$ or smaller. 

Finally, $C_3$ is only of order $n^{2v-3}$ if $H$ contains triangles. If $H$ does not contain triangles, then the only way for two copies of $H$ to have three edges in common is to have four or more vertices in common. Thus, if $H$ does not contain
triangles, then $C_3 = O(n^{2v-4})$ and we are left with 
\be Var(T_H)_{n,p} = C_2p^{2\ell-2}(1-p)^2 + O(n^{2v-4}). \ee
\end{proof}
Simple use of Chebychev's inequality ellucidates the terms of different growth:
\begin{cor}
\be
\frac{T_H(n)}{n^v}\to \frac{1}{|S_H|}p^\ell,\ \ \ \ 
n\Big[{\frac{T_H(n)}{n^v} - \frac{1}{|S_H|}p^\ell}\Big]\to 
\frac{1}{|S_H|}\frac{v(v-1)}{2}p^\ell
\ee
where the random variables are converging in probability.
\end{cor}

We now address the choice we made to use $n$ to measure the `size' of
our random system $G_n$, which was then used when identifying
`surface' effects. The probability distribution on $G_n$ is based on fixing the
number of edges that can appear in the graphs of $G_n$ which we allow,
the graphs which appear in our analysis. In this sense the size of
$G_n$ is perhaps more properly $N_n={n \choose 2}$, as the constraint strictly
limits the fraction of the possible $N_n$ possible edges.
If we rewrite our expansions
of the mean and variance of $T_H$ in powers of $N_n$
we get:
\begin{eqnarray}
\langle T_H \RA &=& \frac{1}{|S_H|}[2^{\frac{v}{2}}N_n^{\frac{v}{2}}- 
2^{\frac{v}{2}-\frac32}v(v-2)N_n^{\frac{v}{2}-\frac{1}{2}}]+ O(N_n^{\frac{v}{2}-1}),\cr 
Var(T_H)_{n,p} & =& O(N_n^{v-\frac{3}{2}}).
\end{eqnarray}
This decomposition of $\langle T_H \rangle_{n,p}$ is somewhat different from that
of equation (\ref{thm1eq}), but the standard deviation of $T_H$ 
still has a growth rate,
$N_n^{\frac{v}{2}-\frac{3}{4}}$, that is smaller than
the subleading term in the expansion of the mean of $T_H$. The precise size
of the surface term depends on the choice of size parameter, but the 
existence of a surface effect is unambiguous. 

We will address this issue again in the next section, 
and again in the Conclusion.

\section{A random graph model with independent edges}\label{SecIndep}

Now we turn to the model defined by having all edges appear independently with
probability $p$. (This model is also often called `Erd\"os-R\'enyi', despite being introduced in \cite{Gil}.)
If one 
identifies edges with coin flips the model can be understood as a coin
flip model in which one
focuses on random variables $T_H$ that are not easily described
in the standard setting of coin flips. This presentation
makes it easy to see how adding dependence to the coin flips, through 
fixing the fraction of heads, affects these `graph theoretic' random variables.

In this model, the total number $T_\E$ of edges is a random variable with 
mean $N_np$ and variance $N_np(1-p)$. However, this model can also be used to mimick the
model of the last section with a sharp constraint on the number of
edges, using a residual variance, as we shall see.
The variable $T_H$ is correlated with $T_\E$,
with correlation coeffiecient
\be r = \frac{Cov(T_H, T_\E)_{n,p}}{\sqrt{Var(T_H)_{n,p}Var(T_\E)_{n,p}}}.
\ee
A common interpretation of $r$ is that a fraction $r^2$ of the variance of 
$T_H$ in the dependent-edge model can be ``explained'' by the correlation with $T_\E$, and that the remaining 
{\em residual variance} of $T_H$ is 
\be ResVar(T_H)_{n,p} = (1-r^2) Var(T_H)_{n,p} = Var(T_H)_{n,p} - 
\frac{Cov(T_H, T_\E)_{n,p}^2}{Var(T_\E)_{n,p}}.
\ee
If we model $T_H$ as a linear function of $T_\E$ plus a residual piece that
is uncorrelated to $T_\E$, then $ResVar(T_H)_{n,p}$ is the variance of this 
residual piece. That is, $ResVar(T_H)_{n,p}$ is the variance we should expect
if we further constrain our system to have a specific value of $T_\E$, as
in the previous section.  

\begin{thm}\label{independenttheorem}
In the independent-edge model, the expectation, variance, and
residual variance of $T_H$ 
are given by:
\begin{eqnarray}\label{thm2eq}
 \langle T_H \rangle_{n,p} &=& c_n p^\ell = \frac{1}{|S_H|}[n^v -\frac{v(v
    -1)}{2}n^{v-1}]p^\ell +O(n^{v-2}) \cr 
Var(T_H)_{n,p} &=& \sum_{k} C_k (p^{2\ell-k} - p^{2\ell})  = O(n^{2v-2})\cr 
ResVar(T_H)_{n,p} & = & \sum_k C_k \left (
p^{2\ell-k}-p^{2\ell} -k(p^{2\ell-1}-p^{2\ell})\right ) \cr 
& = & 
C_2 p^{2\ell-2}(1-p)^2 + C_3 p^{2\ell-3}(1-3p^2+2p^3) + O(n^{2v-4})
\cr &=& O(n^{2v-3}).
\end{eqnarray}

\end{thm}

The independent-edge model gives the same
results for the mean of $T_H$, up to unimportant lower-order
corrections, as the dependent-edge model. However the 
variance is one power of $n$ larger than in the dependent-edge
model, so the subleading term in the expansion of the mean of $T_H$
is of the same order, $O(n^{v-1})$, as the standard deviation of
$T_H$, and we say the independent-edge model does not have a surface term.
Not surprisingly, the {\em residual} variance
of $T_H$ in the independent-edge model matches
the variance of $T_H$ in the dependent-edge model.

\begin{proof}
  The calculation is essentially the same as in the dependent-edge model, only
  with $P(\EN,N_n,k)$ replaced by $p^k$. Since there are $c_n$
  configurations for $H$, each with probability $p^\ell$, the
  expectation of $T_H$ is $c_n p^\ell$. We then have 
  \be \langle T_H \rangle_{n,p}^2 = c_n^2 p^{2\ell} = \sum_{k=0}^\ell C_k p^{2\ell}. \ee 
  As for
  $\langle T_H^2 \rangle_{n,p}$, each of the configurations with $k$
  overlapping edges has probability $p^{2\ell-k}$, so 
  \be \langle T_H^2 \rangle_{n,p} = \sum_{k=0}^\ell C_k p^{2 \ell-k}. \ee 
  Subtracting, we get 
  \be Var(T_H)_{n,p} = \langle T_H^2 \rangle_{n,p} - \langle T_H \rangle_{n,p}^2
  = \sum_{k=1}^\ell C_k (p^{2\ell-k}-p^{2\ell}). \ee
This sum is dominated by the $k=1$ term, which scales as $n^{2v-2}$. 

  To get the covariance of $T_H$ and $T_\E$ we must compute the number
  of ways to have an $H$ and a special edge (representing
  $T_\E$). There are $c_n(N_n-\ell)$ ways to have the edge be disjoint from
  the edges of $H$, and $c_n\ell$ ways to have the special edge be one
  of the edges of $H$. Thus
\begin{eqnarray} \langle T_H T_\E \rangle_{n,p} & = & c_n(N_n-\ell)p^{\ell+1} + c_n\ell p^\ell, \cr 
\langle T_H \rangle_{n,p} \langle T_\E \rangle_{n,p} & = & c_n(N_n-\ell)p^{\ell+1} + c_n \ell p^{\ell+1}, \cr 
Cov(T_H, T_\E)_{n,p} & = & c_n \ell (p^{\ell} - p^{\ell + 1})= c_n \ell p^\ell(1-p). 
\end{eqnarray}
We then have 
\be \frac{Cov(T_H, T_\E)_{n,p}^2}{Var(T_\E)_{n,p}} = \frac{c_n^2 \ell^2}{N_n}p^{2\ell-1}(1-p). \ee
However, $c_n^2 \ell^2/N = \sum_{k=0}^\ell k C_k$, so 
\be Var(T_H)_{n,p} - \frac{(Cov(T_H, T_\E)_{n,p}^2}{Var(T_\E)_{n,p}} = \sum_{k=0}^\ell C_k \left (p^{2\ell-k}-p^{2\ell} - k(p^{2\ell-1}-p^{2\ell}) \right ).
\ee
The $k=0$ and $k=1$ terms are identically zero, the terms with $k>3$ are of order $O(n^{2v-4})$ or smaller, and what is left
is $C_2 p^{2\ell-2}(1-p)^2 + C_3 p^{2\ell-3} (1-3p^2 + 2p^3)$. 
\end{proof}

Note that some of the equations in the theorem do not need lower-order
corrections. The expectation agrees with the dependent-edge model up
to order $O(n^{v-2})$, while the {\em residual} variance of the independent-edge
model agrees with the variance of the dependent-edge model up to order
$O(n^{2v-4})$.

The scale and relative lack of statistical significance of the
surface term is unaffected by the choice of measure of the size of $G_n$. In terms of $N_n$, we have
\begin{eqnarray}
\langle T_H \RA &=& \frac{1}{|S_H|}[2^{\frac{v}{2}}N_n^{\frac{v}{2}}- 2^{\frac{v}{2}-\frac32} v(v-2) N_n^{\frac{v}{2}-\frac{1}{2}}]+ O(N_n^{\frac{v}{2}-1}),\cr 
Var(T_H)_{n,p} & =& O(N_n^{v-1}), \cr 
ResVar(T_H)_{n,p} & = & O(N_n^{v-\frac32}),
\end{eqnarray}
so the standard deviation of $T_H$ has a growth rate, 
$N_n^{\frac{v}{2}-\frac{1}{2}}$, equal to that of the second 
term in the expansion of the mean of $T_H$.

\section{2-stars, triangles and squares}

Now we work out three examples, specifically where $H$ is a graph with
3 vertices and 2 edges (often called a ``2-star'' or a ``cherry''), where
$H$ is a triangle, and where $H$ is a square.  

\subsection{2-stars}

If $H$ is a 2-star, then $c_n= n(n-1)(n-2)/2 = (n^3-3n^2+2n)/2$, $C_2=c_n$, 
and $C_3=0$. Thus the expectation and variance 
in the dependent-edge model (i.e. the first model) are
\begin{eqnarray}
\langle T_H \rangle_{n,p} & = & c_n P(\EN,N_n,2) \cr 
& = & \left ( \frac{n^3}{2} - \frac{3n^2}{2} + O(n) \right )
\left (p^2 + O(n^{-2}) \right ) \cr 
& = & \frac{p^2}{2}n^3 - \frac{3p^2}{2}n^2 + O(n), \cr 
Var(T_H)_{n,p} & = & C_2 p^2(1-p)^2 + O(n^2) \cr 
& = &  \frac{p^2(1-p)^2}{2}n^3 + O(n^2).
\end{eqnarray}

For the independent-edge model, we also need to compute $C_1$, which
works out to equal $2n(n-1)(n-2)(n-3)=2n^4-12n^3+O(n^2)$. 
The variance is then
\begin{eqnarray}
Var(T_H)_{n,p} & = & C_1 (p^5-p^6) + C_2(p^4-p^6) \cr 
& = & (2n^4-12n^3)(p^5-p^6) + \frac{n^3}{2}(p^4-p^6) + O(n^2) \cr 
& = & 2(p^5-p^4)n^4 + \frac{p^4-24p^5+23p^6}{2}n^3 + O(n^2),
\end{eqnarray}
and the residual variance is 
\begin{eqnarray}
ResVar(T_H)_{n,p} & = & C_2 p^2(1-p)^2 \cr 
& = & \frac{p^2(1-p)^2}{2}n^3 + O(n^2).
\end{eqnarray}

\subsection{Triangles}

When $H$ is a triangle, our relevant combinatorial factors are:
\begin{eqnarray}
c_n  =&  \frac{n(n-1)(n-2)}{6}  &= \frac{n^3}{6} - \frac{n^2}{2} + O(n), \cr 
C_1  =& \frac{n(n-1)(n-2)(n-3)}{2} 
 &= \frac{n^4}{2} -3n^3 + O(n^2), \cr 
C_2  && =  0,   \cr 
C_3  =&  c_n & = \frac{n^3}{6} + O(n^2). 
\end{eqnarray}

In the dependent-edge model, we have 
\begin{eqnarray}
\langle T_H \rangle_{n,p} & = & c_n P(\EN,N_n,3) \cr 
& = & \frac{p^3}{6}n^3 - \frac{p^3}{2}n^2 + O(n), \cr 
Var(T_H)_{n,p} & = & C_3 p^3(1-3p^2+2p^3) + O(n^2) \cr 
& = &  \frac{p^3(1-3p^2+2p^3)}{6}n^3 + O(n^2).
\end{eqnarray} 
 
In the independent-edge model we have 
\begin{eqnarray}
\langle T_H \rangle_{n,p} & = & c_n p^3 \cr 
& = & \frac{p^3}{6}n^3 - \frac{p^3}{2}n^2 + O(n), \cr 
Var(T_H)_{n,p} & = & C_1(p^5-p^6) + C_3(p^3-p^6)  \cr 
& = &  \frac{p^5-p^6}{2}n^4 + \frac{p^3-18 p^5 + 17p^6}{6} n^3
+ O(n^2), \cr 
ResVar(T_H)_{n,p}&=& C_3p^3(1-p)^2(1+2p) \cr 
& = & \frac{p^3(1-p)^2(1+2p)}{6}n^3 + O(n^2).
\end{eqnarray}

\subsection{Squares}

If $H$ is a square, then $c_n=n(n-1)(n-2)(n-3)/8$, since we are picking 
4 points and the group of symmetries of the square is the dihedral group
of order 8.  We then compute
\begin{eqnarray}
C_1 & = & \frac{n(n-1)(n-2)(n-3)(n-4)(n-5)}{2} \cr 
& = & \frac{n^6}{2} - \frac{15n^5}{2} + O(n^3), \cr
C_2 & = & \frac{n(n-1)(n-2)(n-3)(n-4)}{2}+ \frac{n(n-1)(n-2)(n-3)}{4} \cr 
& = & \frac{n^5}{2} + O(n^4), \cr 
C_3 & = & 0, \cr 
C_4 & = & c_n = O(n^4).
\end{eqnarray}
The first term in $C_2$ comes from having two consecutive edges shared across
the two squares, while the second comes from sharing non-consecutive edges.

In the dependent-edge model we then have
\begin{eqnarray}
\langle T_H \rangle_{n,p} & = & c_n P(\EN,N_n,4) \cr 
& = & \frac{p^4}{8}n^4 - \frac{3p^4}{4}n^3 + O(n^2), \cr 
Var(T_H)_{n,p} & = & C_2 p^6(1-p)^2 + O(n^4) \cr 
& = &  \frac{p^6(1-p)^2}{2}n^5 + O(n^4). 
\end{eqnarray} 

In the independent-edge model, we have 
\begin{eqnarray}
\langle T_H \rangle_{n,p} & = & c_n p^4 \cr 
& = & \frac{p^4}{8}n^4 - \frac{3p^4}{4}n^3 + O(n^2), \cr 
Var(T_H)_{n,p} & = & C_1(p^7-p^8) + C_2(p^6-p^8)+ C_4(p^4-p^8)  \cr
& = & \frac{p^7-p^8}{2}n^6 + \frac{p^6 -15p^7 + 14p^8}{2}n^5 
+ O(n^4),  \cr
ResVar(T_H)_{n,p}&=& C_2 p^6(1-p)^2 + O(n^4) \cr 
& = & \frac{p^6(1-p)^2}{2}n^5 + O(n^4).
\end{eqnarray}
\section{Block Models}\label{block}

In this section we sketch a more complex version of the previous
models, in which there are vertices of various colors. More specifically,
we consider 
colored graphs on $B$ colors, where the number $n_1, \ldots, n_B$ of vertices
of each color is fixed. We imagine a limit in which all the $n_i$'s go 
to infinity along a fixed line in $\R^B$. 
In the dependent-edge version of this model, we fix
the number $E_{n,ij}$ of edges between vertices of colors $i$ and $j$. In
the independent-edge version of this model, we fix the probability $p_{ij}$
of each such edge. 

In the interest of brevity, we merely sketch the results. (Precise statements and proofs will appear in a subsequent
paper.) 
We have 
\be T_H = \sum_\alpha T_{H_\alpha}, \ee
where $\alpha$ indexes all the possible colorings of $H$. Each $T_{H_\alpha}$
has its expectation and variance described by expansions similar 
to (\ref{thm1eq}) or (\ref{thm2eq}), and similar formulas apply to the
covariances of different $H_\alpha$'s. 
As before, $\langle T_{H_\alpha} \rangle_{n,p}$ always scales as $n^v$, 
while the 
(co)variances of the $T_{H_\alpha}$'s in the dependent-edge model scale as
$n^{2v-3}$, as do the 
residual (co)variances in the independent-edge model. The total (co)variances
in the independent-edge model scale as $n^{2v-2}$. 
As before, the expectations are the
same in the two models (up to $O(n^{v-2})$ corrections), and the residual
variance in the independent-edge model equals the variance in the 
dependent-edge model, up to $O(n^{2v-4})$ corrections. 

The combinatorial factors
$c_n, C_0, C_1,$ etc. are different for different values of $\alpha$, as
are the probabilistic functions that replace $p^\ell$ and $p^{\ell-k}$. 
As a result, $\langle T_H \rangle_{n,p}$ cannot be written as 
a single function of the $n_i$'s times a single function of the $p_{ij}$'s. 
To get an asymptotic understanding of $\langle T_H \rangle_{n,p}$, 
it is necessary
to isolate all the different terms that are bigger than the 
standard deviation. That is, the leading terms of order $n^v$ and the 
surface corrections of order $n^{v-1}$. 

In the dependent-edge model, the subleading terms in the expansion of 
$\langle T_H \rangle_{n,p}$ are $O(\sqrt{n})$ larger than the standard deviation.
Regardless of whether we measure the size of our system in terms of $n$, 
${n \choose 2}$, or some other yardstick, there is an unambiguous surface
effect. 
By contrast, in the independent-edge model 
the subleading terms in the expansion of 
$\langle T_H \rangle_{n,p}$ are of the same order 
as the standard deviation. 

\section{Conclusion}
We considered a sequence $G_n$ of spaces of random graphs 
 through which we study the growth rates of certain random
counts, for instance triangles. The probability distributions on $G_n$
are defined by strongly restricting the count of edges, and
this restriction turns out to reduce the randomness in the counts of
triangles, and indeed any other graph $H$, to such an extent that a
surface phenomenon is produced (Theorem~\eqref{dependenttheorem}):
a lower order constant correction to the mean of the
count of $H$, with growth rate larger than that of the fluctuations. Without the constraint there is no surface effect (Theorem~\eqref{independenttheorem}).

This work was motivated by previous studies of random graph models in
which the randomness is produced by restrictions on the counts of two
or more graphs, say both edges and triangles, and then counts of
other graphs $H$ are studied~\cite{RS1,RS2,RRS1,KRRS1,RRS2,KRRS2,K}.
(When one has two or more count
restrictions they can interfere and produce `phase transitions',
drastic sensitivity in the {\it highest order terms} of counts for $H$, encoded in what
is called the entropy.) In those random graph models the highest order terms in
the counts of graphs $H$ turn out to be easily computable because the 
highest order terms are represented by block models \cite{KRRS2}. This is one of
the reasons we
have included block models in Section \ref{block}. 

Some of that modelling, for instance the edge/triangle model, was
explicitly performed to help understand features in statistical
mechanics. Statistical mechanics was created by Boltzmann and Gibbs
based on two conservation laws, the fact that the sum of the energies
of all the particles, and the sum of the masses of all the particles,
are dynamically conserved and therefore can each be rigorously fixed
as adjustable parameters. The way we produced the probability
distribution on our $G_n$ is an explicit copy of this, but only using
the mass conservation. We would have
liked to restrict two or more graphs (to study phase transitions) but
were not able to control the combinatorics to look for surface
effects when the leading order terms were so sensitive.

What was done here could all be done, in principle, in other combinatorial settings, for instance the sequence of spaces $P_n$ of permutations on $n$ objects. There is some literature~\cite{KKRW} on random permutations in which constraints are put on the counts of two or more `patterns', in order to study interactions between the constraints in the highest order terms in the expansions of the counts of other patterns, i.e.\ phase transitions. It would be of interest to explore the existence of surface effects in random pattern counts using only one pattern restriction.



\begin{thebibliography}{1234}

\bibitem{JLR}  S.\ Janson, T.\ Luczak and A.\ Rucinski, {\em Random Graphs},
John Wiley, New York, 2000

\bibitem{B} B. Bollobas, {\em Random Graphs, 2nd ed.}
Cambridge University Press, Cambridge, 2001.

\bibitem{ER} P.~Erd\H{o}s and A.~R\'enyi, On Random Graphs, 
{\em Publ. Math.} 6 (1959) 290--297.

\bibitem{Gil} E.N.~Gilbert, Random Graphs, {\em Ann. Math. Stat.} 30 (1959) 1141-1144.

\bibitem{KRRS1} R.\  Kenyon, C.\  Radin, K.\  Ren and L.\  Sadun, Multipodal
  structure and phase transitions in large constrained graphs,
  {\em J. Stat. Phys.} 168 (2017) 233-258.

\bibitem{KRRS2} R.\  Kenyon, C.\  Radin, K.\  Ren and L.\  Sadun, Bipodal
  structure in oversaturated random graphs, {\em
  Int. Math. Res. Notices}, 2016 (2016) rnw261.

\bibitem{RS1} C. Radin and L. Sadun, Phase transitions in a complex
  network, {\em J. Phys. A: Math.
Theor.} 46 (2013) 305002.

\bibitem{RS2} C. Radin and L. Sadun, Singularities in the entropy of asymptotically large simple
graphs, {\em J. Stat. Phys.} 158 (2015) 853-865.

\bibitem{RRS1} C. Radin, K. Ren and L. Sadun, The asymptotics of large
  constrained graphs, {\em J.
Phys. A: Math. Theor.} 47 (2014) 175001.

\bibitem{RRS2} C. Radin, K. Ren and L. Sadun, A symmetry breaking transition in the edge/triangle
network model, {\em arXiv:1604.07929v1} (2016).

\bibitem{K} H.\ Koch,
Vertex order in some large constrained random graphs,
{\em SIAM J. Math. Anal.} 48 (2016) 2588-2601.

\bibitem{KKRW} R.\  Kenyon, D. \ Kr\'{a}l', C.\  Radin and P.\  Winkler,
Permutations with fixed pattern densities, {\em arXiv:1506.02340v2.} (2015).

\end{thebibliography}
\end{document}